\documentclass{ifacconf}
\usepackage{amsmath}
\usepackage{amssymb}
\usepackage{graphicx}      % include this line if your document contains figures
\usepackage{natbib}        % required for bibliography
\usepackage{arydshln}
\usepackage{mathtools}
\usepackage[super]{nth}
\usepackage{xfrac}
\usepackage[shortlabels]{enumitem}
\usepackage{xspace}
\usepackage{pifont}

\makeatletter
\def\old@comma{,}
\catcode`\,=13
\def,{%
  \ifmmode%
    \old@comma\discretionary{}{}{}%
  \else%
    \old@comma%
  \fi%
}
\makeatother

\allowdisplaybreaks

\newcommand{\tp}{^{\top}}

\newcommand{\iv}{^{-1}}
\newcommand*\mystrut[1]{\vrule width0pt height0pt depth#1\relax}
\newcommand{\MATLAB}{\textsc{Matlab}\xspace}
\DeclarePairedDelimiter{\norm}{\lVert}{\rVert}

\begin{document}
\begin{frontmatter}

\title{Numerically Reliable Brunovsky Transformations\thanksref{footnoteinfo}} 

\thanks[footnoteinfo]{This project has received funding from the European Union’s 2020 research and innovation programme under the Marie Skłodowska-Curie grant agreement No. 953348 ELO-X. }

\author[First]{Shaohui Yang} 
\author[First]{Colin N. Jones} 

\address[First]{Automatic Control Lab, EPFL, Lausanne, Switzerland  (e-mail: shaohui.yang@epfl.ch, colin.jones@epfl.ch)}

\begin{abstract}                % Abstract of 50--100 words
The Brunovsky canonical form provides sparse structural representations that are beneficial for computational optimal control, yet existing methods fail to compute it reliably. 
We propose a technique that produces Brunovsky transformations with substantially lower construction errors and improved conditioning. 
A controllable linear system is first reduced to the staircase form via an orthogonal similarity transformation. 
We then derive a simple linear parametrization of the transformations yielding the unique Brunovsky form. 
Numerical stability is further enhanced by applying a deadbeat gain before computing system matrix powers and by optimizing the linear parameters to minimize condition numbers.
\end{abstract}


\end{frontmatter}
%===============================================================================

\section{Introduction}
In systems and control theory, the Brunovsky canonical form boosts both theoretical analysis and computational design. 
It forms a cornerstone of geometric control theory and is closely related to feedback linearization in nonlinear control~\citep{slotine1991applied}.
Its independent chain-of-integrator dynamics allow each input-output pair to be controlled individually.
From a state-space perspective, its sparsity pattern accelerates the matrix computation when solving model predictive control (MPC) problems for linear time-invariant (LTI) systems~\citep{yang2025brunovsky}. 
The computational viewpoint has long been ignored, in part due to the lack of a numerically reliable algorithm for transforming any controllable LTI system into its unique Brunovsky form. 
Existing constructive methods produce single transformations that suffer from ill-conditioning. 

In this paper, we bridge this gap by leveraging decoupling theory and key properties of square decouplable systems.
Given a controllable LTI system, we first transform it into staircase form.
We then show that the transformations leading to the Brunovsky form admit a remarkably simple \textit{linear} parametrization.
Finally, we introduce two strategies to enhance numerical stability: applying the deadbeat control gain before computing matrix powers, and optimizing the condition numbers of the transformations with respect to the linear parameters.
Our approach differs from prior work in that we consider a \textit{family}—though not necessarily the entire family—of admissible transformations rather than constructing a single feasible instance.
Moreover, it shows that, within the specific context of controllable systems, a nilpotent Jordan structure can be constructed using a well-conditioned transformation, thereby circumventing the classical ill-posedness that affects general Jordan canonical form computation. 

Notation: $\text{col}(X_i)_{i=1}^m$ denotes stacking $m$ matrices vertically, i.e., 
$\begin{bmatrix}
    X_1\tp & X_2\tp & \dots & X_m\tp 
\end{bmatrix}\tp$.
$\text{blkdiag}(X_i)_{i=1}^m$ denotes placing $m$ matrices on the diagonal, i.e., 
$\begin{bsmallmatrix}
    X_1 & \\
    & X_2 \\
    & & \ddots \\
    & & & X_m
\end{bsmallmatrix}$. 
$[A_{ij}]$ means the matrix $A$ is composed of blocks. 
$I_j$ denotes identity matrix of size $j \times j$.
$0_{j \times k}$ denotes zero matrix of size $j \times k$. 
$[m]$ denotes the set of integers $\{1, 2, \dots, m\}$. 

\section{Roadmap}
\begin{enumerate}\addtocounter{enumi}{-1}
    \item In Section~\ref{sec:prelim}, we review existing methods from the last century on transformations to different special forms of linear controllable systems. 
    \item We define in Sections~\ref{subsec:parametrize-C} and~\ref{subsec:parametrize-Brunovsky} a set of linearly parameterized transformations that all map to the Brunovsky form, and we conjecture that the set contains all such transformations.
    \item Because the set is linearly parameterized, we pose an optimization problem in Section~\ref{subsec:opt-numerical} to minimize the condition number of the resulting transformation. Auxiliary numerical improvements are integrated. 
    \item The experiment in Section~\ref{sec:num-exp} demonstrates that the proposed method achieves numerical reliability superior to existing approaches by several orders of magnitude on a low-dimensional system.
\end{enumerate}

\section{Preliminaries}\label{sec:prelim}

A controllable LTI multivariate system, whether in continuous time or discrete time, is represented by 
\begin{equation}\label{eq:def-linear-system}
    x^+ = Ax + Bu,
\end{equation}
where $x \in \mathbb{R}^{n}$ and $u \in \mathbb{R}^{m}$ denote the state and input vectors, $x^+$ denotes the state derivative or the next state, and $A \in \mathbb{R}^{n \times n}, B \in \mathbb{R}^{n \times m}$ are constant. 
We adopt the common assumptions that the columns of $B$ are linearly independent and $n \geq m$ holds. 

\subsection{Controllability Indices}
The controllability matrix, which is defined by
\begin{equation}
    \mathcal{C} = \begin{bmatrix}
        b_1 & \dots & b_m & Ab_1 & \dots & A b_m & \dots & A^{n-1}b_1 & \dots A^{n-1}b_m 
    \end{bmatrix}, \notag
\end{equation}
has rank $n$, where $b_i$ is the $i$th column of $B$. 
The controllability indices $\{\mu_i\}_{i=1}^{m}$ are defined so that 
\begin{equation}\label{eq:def-modified-ctrbility-matrix}
    \bar{\mathcal{C}} = \begin{bmatrix}
        b_1 & Ab_1 & \dots & A^{\mu_1-1}b_1 & \dots & b_m & \dots & A^{\mu_m-1}b_m 
    \end{bmatrix}
\end{equation}
is an invertible $n\times n$ matrix. The integer $\mu_i$ denotes the number of columns involving $b_i$ in the set of the first $n$ linearly independent columns found in $\mathcal{C}$ when moving from left to right. 
$\sum_{i=1}^{m} \mu_i = n$ holds. 
The set is unique if arranged in descending order $\mu_1 \geq \mu_2 \geq \dots \geq \mu_m \geq 1$. 
The controllability index $\mu$ is the maximum of: $\mu = \underset{i}{\max} \hspace{0.2em} \mu_i = \mu_1$. 
The set $\{\mu_i\}_{i=1}^{m}$ is invariant under invertible state and input transformations $T \in \mathbb{R}^{n \times n}, G \in \mathbb{R}^{m \times m}$, and state feedback $F \in \mathbb{R}^{m \times n}$~\citep{antsaklis2006linear}. 
$T, G$, and $F$ deliver the following change of coordinates
\begin{equation}
        x = T\iv z, \quad
        u = Fx + Gv,
\end{equation}
where $z \in \mathbb{R}^n$ and $v \in \mathbb{R}^m$ denote the state and input vectors in the new coordinate. 
System~\eqref{eq:def-linear-system} becomes
\begin{equation}\label{eq:linear-system-new-coordinate}
    z^+ = T (A+BF) T\iv z + TBG v. 
\end{equation}
Equivalently speaking, the controllability indices of system~\eqref{eq:def-linear-system} are the same as system~\eqref{eq:linear-system-new-coordinate}. 
A numerically stable approach to compute $\{\mu_i\}_{i=1}^{m}$ is to construct an orthogonal state transformation $U \in \mathbb{R}^{n \times n}$ such that
\begin{equation}\label{eq:trans-staircase}
    U A U\tp = A_s, \quad U B = B_s, 
\end{equation}
are in staircase form~\citep{van1979computation}, get a set of integers from the sizes of matrix blocks, and compute the dual partition of the set. 
The staircase form is defined as
\begin{equation}\label{eq:def-staircase-form}
    A_s = \begin{bmatrix}
        A_{11}^s & A_{12}^s & \dots & \dots & A_{1,\mu}^s \\
        A_{21}^s & A_{22}^s & \dots & \dots & A_{2,\mu}^s \\
        0 & A_{32}^s & \dots & \dots & A_{3,\mu}^s \\
        0 & 0 & \dots & A_{\mu,\mu-1}^s & A_{\mu,\mu}^s
    \end{bmatrix}, \quad
    B_s = \begin{bmatrix}
        A_{10}^{s} \\
        0 \\
        \vdots \\\
        0
    \end{bmatrix},     
\end{equation}
where $A_{ij}^s \in \mathbb{R}^{\omega_i \times \omega_j}$ and $A_{i,i-1}^s$ have rank $\omega_{i}$.  
The set of integers $\{\omega_i\}_{i=1}^\mu$ is known as the Weyr characteristics in linear algebra terminology.
It satisfies $m = \omega_1 \geq \omega_2 \geq \dots \geq \omega_\mu \geq 1$ and $\sum_{i=1}^\mu \omega_i = n$. 
Both $\{\omega_i\}_{i=1}^\mu$ and $\{\mu_i\}_{i=1}^{m}$ are integer partitions of $n$ and they are conjugate to one another. 
We avoid referring to~\eqref{eq:def-staircase-form} as ``block upper Hessenberg matrix with Weyr block sizes'', since such a description only characterizes $A_s$ while neglecting $B_s$.

\begin{fact}
    Computed from a controllable linear system pair $(A, B)$, neither the orthogonal transformation $U$ in~\eqref{eq:trans-staircase} nor the staircase form $(A_s, B_s)$ in~\eqref{eq:def-staircase-form} are unique.
\end{fact}

\subsection{Two Kinds of Canonical Forms}
The controllable canonical form $(A_c, B_c)$, first proposed by~\cite{luenberger1967canonical}, is represented by
\begin{align}\label{eq:def-controllable-canonical-form}
    A_c &= [A_{ij}^c], \hspace{0.1em}
    A_{ii}^c = \begin{bmatrix}
        0_{ (\mu_i-1) \times 1}  & I_{\mu_i-1} \\
        \multicolumn{2}{c}{*_{1 \times \mu_i}}
    \end{bmatrix}, \hspace{0.1em}
    A_{ij}^c = \begin{bmatrix}
        0_{(\mu_i -1) \times \mu_j} \\
        *_{1 \times \mu_j}
    \end{bmatrix} \nonumber
    \\
    B_c &= \text{col}(B_i^c)_{i=1}^m, \hspace{0.25em}
    B_i^c = \begin{bmatrix}
        \multicolumn{3}{c}{0_{ (\mu_i-1) \times m}} \\
        0_{1 \times (i-1)} & 1_{1 \times 1} & *_{1 \times (m-i)}
    \end{bmatrix},
\end{align}
where $A_{ij}^c \in \mathbb{R}^{\mu_i \times \mu_j}, B_i^c \in \mathbb{R}^{\mu_i \times m}$, the $1$ in the last row of $B_i^c$ occurs at the $i$th column, and $*$ denotes non-fixed entries. 
All controllable systems~\eqref{eq:def-linear-system} can be transformed to~\eqref{eq:def-controllable-canonical-form} with a non-unique state transformation $T$: 
\begin{equation}\label{eq:trans-controllable-canonical}
    T A T\iv = A_c, \quad TB = B_c. 
\end{equation}
A more concise form, the Brunovsky canonical form $(A_b, B_b)$, was first proposed by~\cite{brunovsky1970classification} and is uniquely determined by the set $\{\mu_i\}_{i=1}^{m}$: 
\begin{align}\label{eq:def-brunovsky-canonical-form}
    A_b &= \text{blkdiag}(A_{i}^b), \quad 
    A_{i}^b = \begin{bsmallmatrix}
        0_{ (\mu_i-1) \times 1}  & I_{\mu_i-1} \\
        0_{1 \times 1} & 0_{ 1 \times (\mu_i-1)}
    \end{bsmallmatrix} \in \mathbb{R}^{\mu_i \times \mu_i}, \nonumber \\
    B_b &= \text{blkdiag}(B_{i}^b), \quad 
    B_i^b = \begin{bmatrix}
        0_{ (\mu_i-1) \times 1}  \\
        1_{1 \times 1} \\
    \end{bmatrix} \in \mathbb{R}^{\mu_i \times 1}. 
\end{align}
All controllable systems~\eqref{eq:def-linear-system} can be transformed to~\eqref{eq:def-brunovsky-canonical-form} with non-unique transformations $T, G$, and $F$: 
\begin{equation}\label{eq:trans-brunovsky-canonical}
    T (A+BF) T\iv = A_b, \quad TBG = B_b. 
\end{equation}
A triplet of $(T, F, G)$ satisfying~\eqref{eq:trans-brunovsky-canonical} will be referred to as a \textit{Brunovsky transformation} from now on. 
In linear algebra terminology, $\{\mu_i\}_{i=1}^{m}$ is known as the Segre characteristics of $A_b$ relative to the eigenvalue $0$. 
We avoid referring to~\eqref{eq:def-brunovsky-canonical-form} as ``Jordan canonical form with zero eigenvalues'', since such a description characterizes $A_b$ while neglecting $B_b$.

\subsection{Transformations to the Brunovsky Form}

In the literature, there are two separate routines to assemble Brunovsky transformations: 
\begin{subequations}
\begin{align}
    (A, B) &\overset{\eqref{eq:trans-staircase}}{\longrightarrow} (A_s, B_s) \longrightarrow (A_b, B_b) \label{eq:routine1} \\
    (A, B) &\overset{\eqref{eq:trans-controllable-canonical}}{\longrightarrow} (A_c, B_c) \longrightarrow (A_b, B_b) \label{eq:routine2}
\end{align}
\end{subequations}
\cite{ford1984simple} adopted~\eqref{eq:routine1}, where the \nth{1} step is numerically stable, but the \nth{2} lacks proof of correctness. 

On the contrary, the \nth{2} step of~\eqref{eq:routine2} is straightforward via eliminating the last rows of $A_{ij}^c$ and $B_i^c$ with $F$ and $G$. 
Such $G$ is an upper triangular matrix with ones on the diagonal~\citep{antsaklis2006linear}. The \nth{1} step has many options. 
\cite{luenberger1967canonical} first proposed to leverage a modified version of the controllability matrix $\bar{\mathcal{C}}$ in~\eqref{eq:def-modified-ctrbility-matrix}. 
Define $\sigma_k = \sum_{i=1}^k \mu_i, k \in [m]$. 
Let $q_k \in \mathbb{R}^{1 \times n}$ denote the $\sigma_k$th row of $\bar{\mathcal{C}}\iv$. 
The following $T_L$ achieves~\eqref{eq:trans-controllable-canonical}. 
\begin{equation}\label{eq:luenberger-canonical-transformation}
    T_L \coloneqq \text{col}(Q_k)_{k=1}^m, \quad 
    Q_k \coloneqq \text{col}(q_k A^i)_{i=0}^{\mu_k-1}. 
\end{equation}
\cite{jordan1973efficient} reduced the time complexity compared with Luenberger's method. 
\cite{datta1977algorithm} provided the minimum number of non-fixed entries in $(A_c, B_c)$. 
However, all methods suffer from solving certain forms of $Xy = z$ with an ill-conditioned left-hand side and matrix power computation $A^\mu$. 
Since the eigenvalues of $A$ may have large magnitude, entries of $A^\mu$ may explode. 

\subsection{Deadbeat Control Gain}\label{subsec:deadbeat-preliminary}

For a discrete-time system~\eqref{eq:def-linear-system}, a deadbeat controller is a state feedback law $u = Kx$ such that the closed-loop system $x^+ = (A+BK)x$ has all eigenvalues at the origin, or equivalently, the matrix $A+BK$ is nilpotent, i.e., $(A+BK)^k = 0$ for some integer $k \leq n$. 
All deadbeat control gains form a set defined as
\begin{equation}\label{eq:def-deadbeat-control}
\begin{split}
    \mathcal{D} \coloneqq & \{ K \mid \exists T \text{ s.t. } T(A+BK)T\iv = A_b \text{ in~\eqref{eq:def-brunovsky-canonical-form} with} \\ 
    & \text{ Jordan block sizes decided by } \{\nu_i\}_{i=1}^m \},
\end{split}
\end{equation}
where the admissible Jordan block sizes $\{\nu_i\}_{i=1}^m$ are constrained by the controllability indices $\{\mu_i\}_{i=1}^m$ of the system~\citep{funahashi1992explicit}
\begin{equation}\label{eq:feasible-jordan-block-size}
    \sum_{i=1}^{j} \nu_i \geq \sum_{i=1}^{j} \mu_i, \forall j \in [m]. 
\end{equation}
The following subset of $\mathcal{D}$ has received the most attention
\begin{align}
    \mathcal{D}_{\text{cid}}  \coloneqq & \{ K \mid \exists T \text{ s.t. } T(A+BK)T\iv = A_b \text{ in~\eqref{eq:def-brunovsky-canonical-form} with} \notag \\ 
    & \text{ Jordan block sizes decided by } \{\mu_i\}_{i=1}^m \} \label{eq:def-jordan-block-size-same-MTDC},
\end{align}
because the Jordan block sizes are identical to the controllability indices. 
Starting from the staircase form~\eqref{eq:def-staircase-form}, \cite{van1984deadbeat} constructed a minimum norm solution within $\mathcal{D}_{\text{cid}}$ with an orthogonal state transformation as byproduct. 
\cite{sugimoto1993direct} synthesized one element of $\mathcal{D}_{\text{cid}}$ with no byproduct in a simpler fashion. 
\cite{amin1988parameterization} proved that the entire set $\mathcal{D}_{\text{cid}}$ can be characterized with the least number of parameters linearly, which motivates this paper. 
\cite{yamada1993straightforward} further parametrized the entire set $\mathcal{D}$ covering all feasible $\{\nu_i\}_{i=1}^m$. 
In this paper, we focus on $\mathcal{D}_{\text{cid}}$ rather than the larger $\mathcal{D}$ because the latter does not provide a feedback transformation targeting the Brunovsky form. 

% \subsubsection{Example on the importance of $\mathcal{D}_{\text{cid}}$ over $\mathcal{D}$} Consider a linear system $(A, B)$ in Brunovsky form with $n = 4, m =2, \mu_1 = \mu_2 = 2$. A deadbeat gain $K = \begin{bsmallmatrix}
%     0 & 0 & 0 & 0 \\
%     0 & a & 0 & 0 
% \end{bsmallmatrix}, a \neq 0$ makes $A + BK = \begin{bsmallmatrix}
%     0 & 1 & 0 & 0 \\
%     0 & 0 & 0 & 0 \\
%     0 & 0 & 0 & 1 \\
%     0 & a & 0 & 0 \\
% \end{bsmallmatrix}$. 
% Since $(A+BK)^3 = 0$, it must be similar to a Jordan canonical form with block sizes $(3,1)$ instead of the controllability indices $(2,2)$. 
% All possible state transformations that make $T(A+BK)T\iv = \begin{bsmallmatrix}
%     0 & 1 & 0 & 0 \\
%     0 & 0 & 1 & 0 \\
%     0 & 0 & 0 & 0 \\
%     0 & 0 & 0 & 0 \\
% \end{bsmallmatrix}$ are of the form $T = \begin{bsmallmatrix}
%     0 & 0 & \sfrac{1}{a} & 0 \\
%     0 & 0 & 0 & \sfrac{1}{a} \\
%     0 & 1 & 0 & 0 \\
%     -a & 0 & 0 & 1 \\
% \end{bsmallmatrix}$. The same $T$ makes $TB = \begin{bsmallmatrix}
%     0 & 0 \\
%     0 & \sfrac{1}{a} \\
%     1 & 0 \\
%     0 & 1
% \end{bsmallmatrix}$. However, it is impossible to find an invertible input transformation $G$ such that $TBG = \begin{bsmallmatrix}
%     0 & 0 \\
%     0 & 0 \\
%     1 & 0 \\
%     0 & 1
% \end{bsmallmatrix}$ or $\begin{bsmallmatrix}
%     0 & 0 \\
%     1 & 0 \\
%     0 & 0 \\
%     0 & 1
% \end{bsmallmatrix}$. 
% This counterexample illustrates that setting Jordan block sizes of $A$ that differ from the controllability indices violates the definition of the Brunovsky form and may cause a sparsity mismatch. 

\subsection{Decouplable Systems}

It is common that the linear system~\eqref{eq:def-linear-system} has observables $y = Cx$, where $y \in \mathbb{R}^m$ is the output vector and $C \in \mathbb{R}^{m \times n}$ is constant with full rank $m$. 
Considerable effort has been devoted to designing transformations that, under decoupling theory, ensure that the $i$th output $y_i$ is solely influenced by the $i$th input $u_i$. 
Let $c_i$ denote the $i$th row of $C$. 
The decoupling indices $\{\alpha_i\}_{i=1}^m$ are defined from 
\begin{equation}\label{eq:def-decouple-indices}
    \alpha_i = \text{min} \{ k \mid c_i A^{k-1} B \neq 0, k \in [n] \}. 
\end{equation}
The decoupling matrix $D \in \mathbb{R}^{m \times m}$ is defined by
\begin{equation}\label{eq:def-D-decouple}
    D = \text{col}(d_i)_{i=1}^m, \quad d_i = c_i A^{\alpha_i-1}B. 
\end{equation}
The system $(A, B, C)$ is decouplable if its decoupling matrix $D$ is nonsingular. 
Define the matrix $C^* \in \mathbb{R}^{m \times n}$
\begin{equation}\label{eq:def-C*-decouple}
    C^* = \text{col}(c_i^*)_{i=1}^m, \quad c_i^* = c_i A^{\alpha_i}. 
\end{equation}
Construct matrix $W \in \mathbb{R}^{\alpha \times n}$ with $\alpha = \sum_{i=1}^m \alpha_i \leq n$ by
\begin{equation}\label{eq:def-W-decouple}
    W \coloneqq \text{col}(W_k)_{k=1}^m, \quad 
    W_k \coloneqq \text{col}(c_k A^i)_{i=0}^{\alpha_k-1}. 
\end{equation}
If the system $(A, B, C)$ is decouplable, then applying input transformation $G = D\iv$ with $D$ defined by~\eqref{eq:def-D-decouple}, feedback transformation $F = -D\iv C^*$ with $C^*$ defined by~\eqref{eq:def-C*-decouple}, and state transformation $T = \begin{bsmallmatrix}
    W \\ L
\end{bsmallmatrix}$ with $W$ defined by~\eqref{eq:def-W-decouple} and $L \in \mathbb{R}^{(n-\alpha) \times n}$ ensuring invertibility of $T$ will result in a canonically decoupled system $\big( A_d=T(A-BD\iv C^*) T\iv, B_d= TBD\iv, C_d = C T\iv \big)$~\citep{gilbert1969decoupling}, where 
\begin{align}\label{eq:def-canonically-decoupled-system}
    A_d &= \begin{bsmallmatrix}
        A_1^d & \\
        & A_2^d \\ 
        & & \ddots \\
        & & & A_m^d \\
        A_{r1}^d & A_{r2}^d & \dots & A_{rm}^d & A_{rr}^d
    \end{bsmallmatrix}, \quad
    \begin{aligned}
        r &= m + 1, \\
        A_{rr}^d &\in \mathbb{R}^{(n-\alpha) \times (n-\alpha)}, \\
        A_i^d &= \begin{bsmallmatrix}
        0_{ (\alpha_i-1) \times 1}  & I_{\alpha_i-1} \\
        0_{1 \times 1} & 0_{ 1 \times (\alpha_i-1)}
        \end{bsmallmatrix} \\
        A_{ri}^d &\in \mathbb{R}^{(n-\alpha) \times \alpha_i}, i \in [m],
    \end{aligned} \notag \\
    B_d &= \begin{bmatrix}
        B_1^d & \\
        & B_2^d \\ 
        & & \ddots \\
        & & & B_m^d \\
        B_{r1}^d & B_{r2}^d & \dots & B_{rm}^d
    \end{bmatrix}, \quad
    \begin{aligned}
        B_i^d &= \begin{bmatrix}
        0_{ (\alpha_i-1) \times 1}  \\
        1_{1 \times 1} \\
        \end{bmatrix}, \\
        B_{ri}^d &\in \mathbb{R}^{(n-\alpha) \times 1}, \\
        i &\in [m],
    \end{aligned} \\
    C_d &= \begin{bmatrix}
        C_1^d & & & & 0 \\
        & C_2^d & & & 0 \\ 
        & & \ddots & & \vdots \\
        & & & C_m^d & 0 \\
    \end{bmatrix}, \quad 
    \begin{aligned}
        C_i^d &= \begin{bmatrix}
            1_{1 \times 1} & 0_{(\alpha_i-1)\times 1}
        \end{bmatrix}, \\
        i &\in [m].        
    \end{aligned} \notag
\end{align}

An interesting discovery by~\cite{amin1988parameterization} reveals that if the decoupling indices are identical to the controllability indices, then $F = -D\iv C^* \in \mathcal{D}_{\text{cid}}$. 
The same fact also results in the degeneration of $(A_d, B_d)$ in~\eqref{eq:def-canonically-decoupled-system} to the Brunovsky form in~\eqref{eq:def-brunovsky-canonical-form} and $C_d$ to $\text{blkdiag}(C_i^d)_{i=1}^m$, which forms the basis of our main results.

\section{Main results}

Given a controllable linear system, our objective is to construct transformations that yield the unique Brunovsky form.
Although no observation matrix is present in our setting, it is useful—by analogy with decoupling theory, which assumes one—to reason in the reverse direction.
Starting from a general dense system $(A, B)$ is challenging, so, following many existing works, we begin with the staircase form $(A_s, B_s)$, whose block shapes are exposed and therefore offer a structural entry point.
We then parameterize all matrices $C$ for which the triplet $(A_s, B_s, C)$ is decouplable with preferred indices. 
Finally, we construct Brunovsky transformations that depend linearly on these parameters and optimize over them to obtain an instance with favorable condition numbers.

\subsection{Decouplable Observation Matrix Parametrization}\label{subsec:parametrize-C}

\begin{thm}\label{thm:observation-matrix-parametric-form}
    Given $(A_s, B_s)$ from~\eqref{eq:def-staircase-form}, the set 
    \begin{align}\label{eq:def-observation+decoupling}
        \mathcal{O} = \{ &C \in \mathbb{R}^{m \times n} \mid (A_s, B_s, C) \text{ is decouplable and } \notag \\
        &\forall i \in [m], \alpha_i = \mu_i \}. 
    \end{align}
    is given by the compact parametric form~\eqref{eq:observation+decoupling-construction} constructed through steps (a-c). 
\end{thm}

$\mathcal{O}$ is of vital importance because if the triplet $(A_s, B_s, C)$ has identical decoupling and controllability indices, then the Brunovsky transformations can be constructed via the decoupling theory.
The constructive parametrization of $\mathcal{O}$ is outlined below, followed by a necessity and sufficiency proof. 
The process starts with staircase form $(A_s, B_s)$.
\begin{enumerate}[(a)]
    \item Append $\omega_0 = m$ and $\omega_{\mu+1} = 0$ to the beginning and end of the Weyr characteristics $\{\omega_i\}_{i=1}^m$. 
    \item Let $\epsilon_i = \omega_i - \omega_{i+1} > 0$ for $i= k_1, k_2, \dots, k_g$.
    Start with $k_1 = \mu$ and continue so that $k_j > k_{j+1}$.
    By construction, $\{ k_j \}_{j=1}^g$ contains the $g$ unique elements of $\{\mu_i\}_{i=1}^{m}$ and there are $\epsilon_{k_j}$ instances of $\mu = k_j$ among the controllability indices. 
    $\sum_{j=1}^g \epsilon_{k_j} = m$ thus holds. 
    \item Use the parametric matrices $S^r$ and $S^f$ to construct the following observation matrix: 
    \begin{align}\label{eq:observation+decoupling-construction}
        C &= \text{col}(C_j)_{j=1}^g, \\
        C_j &= \begin{bmatrix}
            \smash{\underbrace{0}_{\omega_1}} & \smash{\underbrace{0}_{\omega_2}} & \dots & \smash{\underbrace{0}_{\omega_{k_j-1}}} & \mystrut{6ex}\smash{\underbrace{S_{k_j}^r}_{\omega_{k_j}}}  & \smash{\underbrace{S^f_{k_j}}_{\sum_{i>k_j}^{m} \omega_i}}
        \end{bmatrix} \in \mathbb{R}^{\epsilon_{k_j} \times n}, \notag
    \end{align}
    where the following $g$ square matrices
    \begin{equation}\label{eq:rank-constraints}
    S_{k_1}^r, \begin{bsmallmatrix}
        A_{k_j+1, k_j}^s \\
        S_{k_j}^r
    \end{bsmallmatrix}, j=2, 3, \dots, g \text{ are invertible}
    \end{equation}
    and $S^f_{k_j}$ are entirely free. $A_{k_j+1, k_j}^s$ are either the subdiagonal blocks of $A_s$ or the top block of $B_s$. 
    The zero blocks of $C_j$ ensure that $\forall i, \alpha_i = \mu_i$. 
    The invertible matrices in~\eqref{eq:rank-constraints} ensure that the decoupling matrix is nonsingular. 
\end{enumerate}
The total degrees of freedom in~\eqref{eq:rank-constraints} is $N_{\mathcal{R}} + N_{\mathcal{F}}$, where
\begin{equation}\label{eq:def-parameter-dof}
    \begin{aligned}
        N_{\mathcal{R}} &= \sum_{j=1}^g \epsilon_{k_j} \omega_{k_j} \text{ depicts the rank-constrained } S^r_{k_j}, \\
        N_{\mathcal{F}} &= \sum_{j=1}^g \epsilon_{k_j} \cdot \sum_{i>k_j}^{m} \omega_i \text{ depicts the free } S^f_{k_j}.
    \end{aligned}
\end{equation}

Before proving the observation matrix parametrization theorem, an auxiliary lemma on the invertibility of the decoupling matrix $D$ is presented. 
\begin{lem}\label{lem:auxiliary}
    Let matrices $A_j, S_j$ be given and satisfy:
    \begin{enumerate}
        \item For $j \in [p]$, $A_j \in \mathbb{R}^{r_j\times r_{j-1}}$ have full row rank $r_j$. $r_0 = r_1 = m$ and $r_j \leq r_{j-1}$. 
        \item For $j \in [p]$, $S_j \in \mathbb{R}^{k_j \times r_j}$. The block row sizes $k_j$ satisfy the telescopic relations: 
        \begin{equation*}
            k_i + r_{i+1} = r_i, i \in [p-1], \quad
            k_p = r_p. 
        \end{equation*}
    \end{enumerate}
    Then the matrix 
    \begin{equation}\label{eq:D-shape-auxiliary-lemma}
        D = \begin{bmatrix}
            S_p A_p A_{p-1} \dots A_2 A_1 \\
            S_{p-1} A_{p-1} \dots A_2 A_1 \\
            \vdots \\
            S_2 A_2 A_1 \\
            S_1 A_1 
        \end{bmatrix} \in \mathbb{R}^{m \times m},
    \end{equation}
    is invertible if and only if the following $p$ square matrices are invertible: 
    \begin{equation}\label{eq:rank-auxiliary-lemma-constraints}
        S_p \text{ and } \begin{bmatrix}
            A_i \\
            S_{i-1}
        \end{bmatrix} \in \mathbb{R}^{r_{i-1} \times r_{i-1}}, i = p, p-1, \dots, 2
    \end{equation}
\end{lem}

\begin{pf}
    See appendix. \hfill $\blacksquare$
\end{pf}

Now we are ready to prove Theorem~\ref{thm:observation-matrix-parametric-form}:
\begin{pf}
    The first half of the proof is to show that if $C$ is constructed via steps (a-c), then it belongs to $\mathcal{O}$. 
    We first show that the decoupling indices of $C$ are the same as $\{\mu_i\}_{i-1}^m$, and then show that the decoupling matrix is invertible. 
    All $\epsilon_{k_j}$ rows of $C_j$ have decoupling indices $k_j$ because the staircase shape of $(A_s, B_s)$ and the zero blocks' positions of $C_j$ guarantee that $\forall j \in [g]$, 
    \begin{itemize}
        \item $C_j A_s^i B_s = 0, i=0, 1, \dots, k_j-2$, and
        \item $C_j A_s^{k_j-1} B_s = S^r_{k_j} \prod_{i=1}^{k_j} A_{k_j+1-i, k_j-i}^s \coloneqq D_j \neq 0$ (product of full rank matrices are nonzero). $S^f_{k_j}$ does not influence $D_j$ because the last $\sum_{i>k_j}^{m} \omega_i$ rows of $A_s^{k_j-1}B_s$ are zero.  
    \end{itemize}
    Since $k_j = \mu_i$ for some $i$ and there are $\epsilon_{k_j}$ instances of such $\mu_i$, the decoupling indices and the controllability indices coincide. 
    The decoupling matrix can be compactly denoted as $D = \text{col}(D_j)_{j=1}^g$, which is nonsingular as a result of the rank constraint~\eqref{eq:rank-constraints} and Lemma~\ref{lem:auxiliary}. \\
    The second half of the proof is to show that if $\tilde{C} \in \mathcal{O}$, then it satisfies~\eqref{eq:observation+decoupling-construction} and~\eqref{eq:rank-constraints}. 
    The sets $\{ k_j \}_{j=1}^g$ and $\{ \epsilon_{k_j} \}_{j=1}^g$ from steps (a-b) will be used directly since they originate from the intrinsic controllability indices. 
    By definition of $\mathcal{O}$, the triplet $(A_s, B_s, \tilde{C})$ has decoupling indices $\alpha_i = \mu_i, \forall i \in [m]$ and its decoupling matrix $\tilde{D}$ is nonsingular. 
    Denote the top $\epsilon_{k_1}$ rows of $\tilde{C}$ as $\tilde{C}_1$. Denote the later $\epsilon_{k_2}$ rows as $\tilde{C}_2$ and so on, until the last $\epsilon_{k_g}$ rows as $\tilde{C}_g$. 
    Assume without loss of generality that 
    \begin{equation*}
        \tilde{C}_j = \begin{bmatrix}
            \smash{\underbrace{X_{j,1}}_{\omega_1}} & \smash{\underbrace{X_{j,2}}_{\omega_2}} & \dots & \smash{\underbrace{X_{j,k_j-1}}_{\omega_{k_j-1}}} & \mystrut{5ex}\smash{\underbrace{\tilde{S}_{k_j}^r}_{\omega_{k_j}}}  & \smash{\underbrace{\tilde{S}^f_{k_j}}_{\sum_{i>k_j}^{m} \omega_i}}
        \end{bmatrix} \in \mathbb{R}^{\epsilon_{k_j} \times n},
    \end{equation*}
    where matrices $X, \tilde{S}$ are of proper sizes. According to~\eqref{eq:def-decouple-indices}\eqref{eq:def-D-decouple}, $\forall j \in [g]$,
    \begin{itemize}
        \item $\tilde{C}_j A_s^i B_s = 0 = \sum_{l=1}^i X_{j,l} M_{l,i} + X_{j,i+1} \prod_{l=0}^i A_{i+1-l, i-l}^s, i=0, 1, \dots, k_j-2$. $M_{l,i}$ is some matrix products of $A_{ij}^s$. 
        When $i=0$, $X_{j,1} A_{10}^s = 0 \Rightarrow X_{j,1} = 0$ because $A_{10}^s$ is of full row rank. 
        When $i=1$, $X_{j,1}M_{1,1} + X_{j,2} A_{21}^s A_{10}^s = 0 \Rightarrow X_{j,2} = 0$ because $X_{j,1} = 0$ and $A_{21}^s A_{10}^s$ is of full row rank. 
        The same logic continues with incremental $i$, so $X_{j,i+1} = 0, i =0, 1, \dots, k_j -2$ because $\prod_{l=0}^i A_{i+1-l, i-l}^s$ has full row rank (specialty of the staircase form), so the sparsity pattern of $\tilde{C}_j$ must be the same as~\eqref{eq:observation+decoupling-construction}. 
        \item As a result of the positions of the zero blocks of $\tilde{C}_j$, $\tilde{C}_j A_s^{k_j-1} B_s = \tilde{S}^r_{k_j} \prod_{i=1}^{k_j} A_{k_j+1-i, k_j-i}^s \coloneqq \tilde{D}_j \neq 0$. 
    \end{itemize}
    The decoupling matrix $\tilde{D} = \text{col}(\tilde{D}_j)_{j=1}^g$ has a shape of~\eqref{eq:D-shape-auxiliary-lemma}. According to Lemma~\ref{lem:auxiliary}, it is nonsingular only if~\eqref{eq:rank-auxiliary-lemma-constraints} (which is equivalent to~\eqref{eq:rank-constraints}) is satisfied. 
    \hfill$\blacksquare$
\end{pf}

\subsubsection{Example of Theorem~\ref{thm:observation-matrix-parametric-form}} 
\begin{figure}
    \centering
    \includegraphics[width=0.5\linewidth]{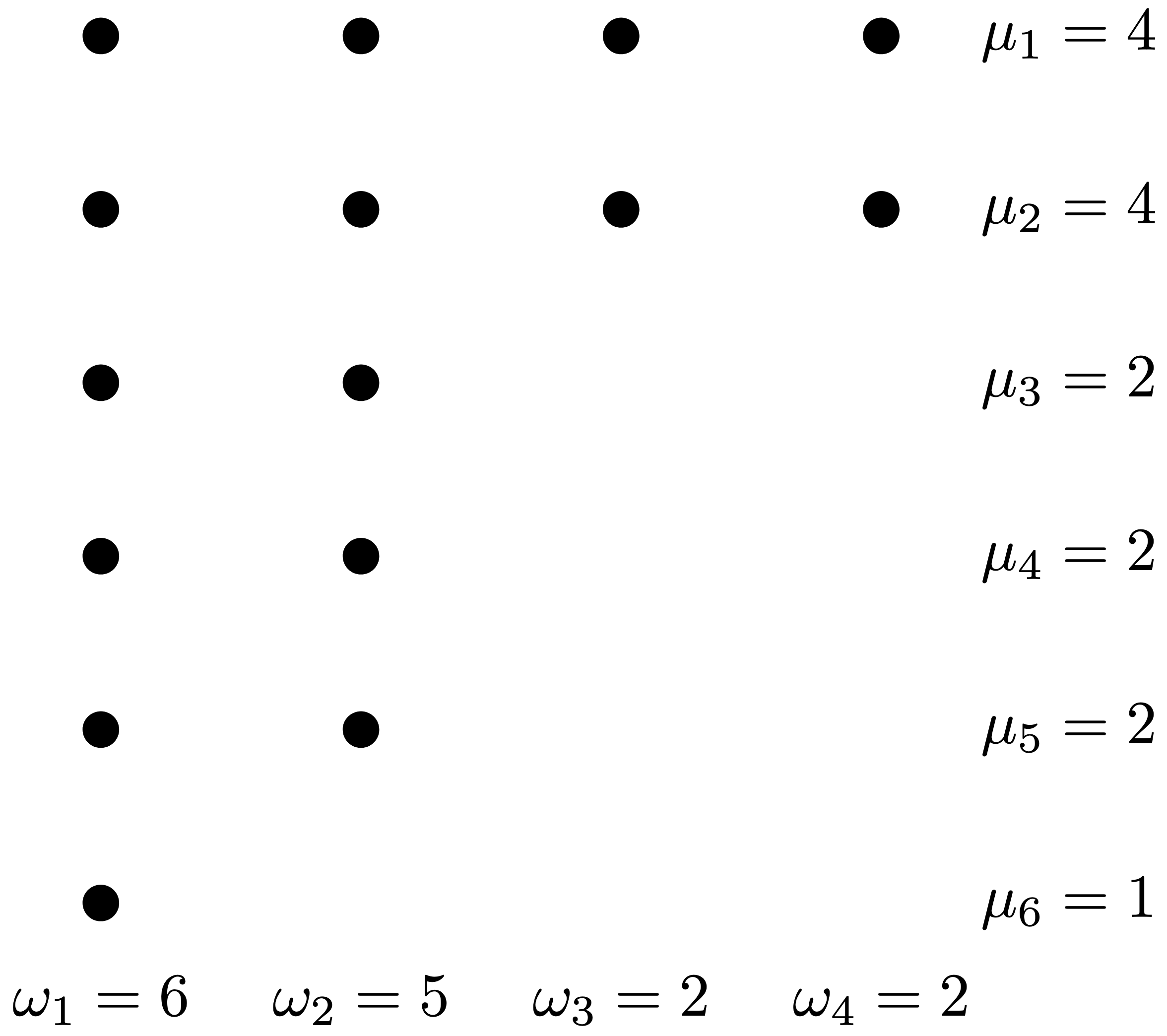}
    \caption{Ferrers diagram. Controllability indices (number of dots per row) and Weyr characteristics (number of dots per column) are conjugate to one another. }
    \label{fig:Ferrers}
\end{figure}

Assume that a linear system has $n = 15$ states and $m=6$ inputs. 
Fig.~\ref{fig:Ferrers} depicts one realization of the controllability indices and the Weyr characteristics. 
Append $\omega_0 = m = 6$ and $\omega_{\mu+1} = \omega_{5} = 0$ to the set of $\omega$. 
$k_1 = \mu = 4$, so $\epsilon_{k_1} = \omega_4 - \omega_5 = 2$. 
Counting from the end, the next $\epsilon_i = \omega_i - \omega_{i+1} > 0$ happens at $i=2$, so $k_2 = 2$ and $\epsilon_{k_2} = \omega_{2} - \omega_{3} = 3$. 
Similarly, $k_3 = 1$ and $\epsilon_{k_3} = \omega_{1} - \omega_{2} = 1$.
Note that the set $\{ k_1=4, k_2=2, k_3=1 \}$ contains the three unique elements of $\{\mu_1=\mu_2 = 4, \mu_3=\mu_4=\mu_5=2, \mu_6 =1\}$ and the set $\{ \epsilon_{k_1}=2, \epsilon_{k_2}=3, \epsilon_{k_3}=1 \}$ reveals that $\mu_i = 4$ happens twice,  $\mu_i = 2$ happens three times, and $\mu_i = 1$ happens once. 
According to Theorem~\ref{thm:observation-matrix-parametric-form}, the following matrix $C$ covers the entire set $\mathcal{O}$:
\begin{equation*}
\begin{aligned}
    C &= \text{col}(C_j)_{j=1}^g \\
    C_1 &= \begin{bmatrix}
        0_{6 \times 2} & 0_{5 \times 2} & 0_{2 \times 2} & S_{k_1}^r  
    \end{bmatrix} & S_{k_1}^r \in \mathbb{R}^{ (\omega_{k_1}=2) \times (\epsilon_{k_1}=2)} \\
    C_2 &= \begin{bmatrix}
        0_{6 \times 3} & S_{k_2}^r & S_{k_2}^f
    \end{bmatrix} 
    & S_{k_2}^r \in \mathbb{R}^{(\omega_{k_2}=5) \times (\epsilon_{k_2}=3)} \\
    C_3 &= \begin{bmatrix}
        S_{k_3}^r & S_{k_3}^f
    \end{bmatrix} & S_{k_3}^r \in \mathbb{R}^{(\omega_{k_3}=6) \times (\epsilon_{k_3}=1)} \\
    S_{k_2}^f &\in \mathbb{R}^{\big( \sum_{i>k_2}^4 = 4 \big) \times (\epsilon_{k_2}=3)} 
    & S_{k_3}^f \in \mathbb{R}^{ \big( \sum_{i>k_1}^4 = 9 \big) \times (\epsilon_{k_3}=1)}
\end{aligned}
\end{equation*}
The decoupling matrix $D = \begin{bsmallmatrix}
        S_{k_1}^r A_{43}^s A_{32}^s A_{21}^s A_{10}^s \\
        S_{k_2}^r A_{21}^s A_{10}^s \\
        S_{k_3}^r A_{10}^s
    \end{bsmallmatrix}$
is nonsingular because the following $g=3$ square matrices are invertible by design: 
$    S_{k_1}^r = S_4^r, \begin{bsmallmatrix}
        A_{32}^s \\
        S_{k_2}^r = S_2^r
    \end{bsmallmatrix}, \begin{bsmallmatrix}
        A_{21}^s \\
        S_{k_3}^r = S_1^r
    \end{bsmallmatrix}$.
The rank-constrained degrees of freedom are $N_{\mathcal{R}} = 25$, obtained as the element-wise product of the \nth{3} and \nth{4} rows of Table~\ref{tab:Ferrers}.
The ``free'' degrees of freedom are $N_{\mathcal{F}} = 21$, computed analogously from the element-wise product of the \nth{5} and \nth{4} rows of Table~\ref{tab:Ferrers}.
\begin{table}[hb]
    \centering
    \begin{tabular}{|c|c|c|c|}
        \hline
        $j$ & 1 & 2 & 3 \\ \hline
        $k_j$ & 4 & 2 & 1 \\ \hline
        $\omega_{k_j}$ & 2 & 5 & 6 \\ \hline
        $\epsilon_{k_j}$ & 2 & 3 & 1 \\ \hline 
        $\sum_{i>k_j}^m \omega_i$ & 0 & 4 & 9 \\ 
        \hline
    \end{tabular}
    \caption{Integer sequences relating to Fig.~\ref{fig:Ferrers}. }
    \vspace{-2em}
    \label{tab:Ferrers}
\end{table}

\subsection{Brunovsky Transformation Parametrization}\label{subsec:parametrize-Brunovsky}

The parametrization of observation matrices $C$ for which $(A_s, B_s, C)$ is decouplable with indices matching the controllability indices provides the foundation for parametrizing Brunovsky transformations, as established in the following theorem.
\begin{thm}\label{thm:parametric-Brunovsky-transformation}
    A linear system in staircase form $(A_s, B_s)$~\eqref{eq:def-staircase-form} can be transformed to its unique Brunovsky form with the parametric state transformation $T$, input transformation $G$, and feedback transformation $F$ defined by
    \begin{align}\label{eq:def-parametric-Brunovsky-transformation}
        T &= \text{col}(T_j)_{j=1}^g & T_j & = \text{col}(T_{j_l})_{l=1}^{\epsilon_{k_j}} \quad T_{j_l} = \text{col}(C_{j_l} A_s^{i})_{i=0}^{k_j-1} \notag \\
        G &= D\iv & D & = \text{col}(C_j A_s^{k_j-1}B_s)_{j=1}^g \\
        F &= -D\iv C^* & C^* & = \text{col}(C_j A_s^{k_j})_{j=1}^g \notag
    \end{align}
    where the parametric $C_j$ satisfies~\eqref{eq:observation+decoupling-construction}, $C_{j_l}$ denotes the $l$th row of $C_j$, and $\{ k_j, \epsilon_{k_j} \}_{j=1}^g$ are constructed via steps (a-b). 
\end{thm}

\begin{pf}
    The decoupling matrix $D$ in~\eqref{eq:def-parametric-Brunovsky-transformation} follows from~\eqref{eq:def-D-decouple}. 
    $C^*$ follows from~\eqref{eq:def-C*-decouple}. 
    $T$ follows from~\eqref{eq:def-W-decouple}. 
    $T, G$, and $F$ achieve $T (A_s+B_sF) T\iv = A_b, TB_sG = B_b$ by the theory of decoupling~\citep{gilbert1969decoupling} because the decoupling indices are identical to the controllability indices, so~\eqref{eq:def-canonically-decoupled-system} degenerates to the Brunovsky form~\eqref{eq:def-brunovsky-canonical-form}. 
    \hfill$\blacksquare$
\end{pf}

\begin{conj}\label{conj:}
    We conjecture that the converse of Theorem~\ref{thm:parametric-Brunovsky-transformation} also holds; namely, if a Brunovsky transformation $(T, F, G)$ transforms the staircase linear system into the Brunovsky form, then there exists an observation matrix $C$ from~\eqref{eq:rank-constraints} that defines the triplet as~\eqref{eq:def-parametric-Brunovsky-transformation}.
    The conjecture is shown graphically by the question mark in~\eqref{eq:conjecture}. 
\end{conj}
\begin{gather}
    \text{Brunovsky transformation } (T, F, G) \text{ for } (A_s, B_s) \notag \\* 
    ? \Downarrow \qquad \Uparrow \checkmark \tag{Conjecture} \label{eq:conjecture} \\*
    \text{Decouplable system } (A_s, B_s, C \in \mathcal{O}) \notag
\end{gather}

Some immediate consequences of Theorem~\ref{thm:parametric-Brunovsky-transformation} and~\eqref{eq:def-parametric-Brunovsky-transformation} are
\begin{enumerate}[1.]
    \item All transformation matrices $(T, F, G)$ are \textit{linear} with respect to the parameters $S^f$ and $S^r$. 
    \item As shown in the proof of Theorem~\ref{thm:observation-matrix-parametric-form}, $D$ and $G$ are solely dependent on $S^r$. 
    \item $F$ and $T$ are dependent on $S^f$ and $S^r$, but $F(S^r, S^f) \in \mathcal{D}_{\text{cid}}$ is overly parametrized. 
    \cite{amin1988parameterization} proved that given a numerical instance of $C \in \mathcal{O}$, $F \in \mathcal{D}_{\text{cid}}$ if and only if $F = -D\iv (C^* + P W)$ where $D, C^*, W$ are defined in the decoupling theory section and $P$ is a sparse parametric matrix with $N_{\mathcal{F}}$ free entries. 
    Their formula for $N_{\mathcal{F}} = nm - \mu_1 - 3 \mu_2 - 5 \mu_3 - \dots - (2m-1)\mu_m$ is equivalent to ~\eqref{eq:def-parameter-dof}. 
    \item A special case of~\eqref{eq:def-parameter-dof} is when all controllability indices are the same, i.e., $\mu_i = \frac{n}{m}, \forall i$. 
    Then $N_{\mathcal{R}} = m^2$ and $N_{\mathcal{F}} = 0$, so $\mathcal{D}_{\text{cid}}$ has an unique element. The observation matrix parametrization in~\eqref{eq:observation+decoupling-construction} degrades to $C = \begin{bmatrix}
            \mystrut{3ex}{\smash{\underbrace{0}_{n - m}}} & S^r
        \end{bmatrix} \in \mathbb{R}^{m \times n}$. Single-Input Single-Output (SISO) systems are obviously such a special case because there is only one controllability index $\mu_1 = \mu = m$. 
\end{enumerate}

\subsection{Numerical Stability Enhancements}\label{subsec:opt-numerical}

Notably, the matrix powers $A_s^i$ are required in~\eqref{eq:def-parametric-Brunovsky-transformation}.
Since $A_s$ is obtained from a unitary similarity transformation of $A$ and may have eigenvalues of arbitrary magnitude, the eigenvalues of $A_s^i$ can grow uncontrollably.
A simple yet effective remedy is to insert a pre-processing step: 
\begin{equation}\label{eq:deadbeat-for-staircase}
    (A_s, B_s) \rightarrow (A_{\tilde{s}} = A_s + B_s K_s, B_s),   
\end{equation}
where $K_s$ is the deadbeat control gain computed from the staircase form~\citep{sugimoto1993direct}.
Since $A_{\tilde{s}}^{\mu} = 0$, the entries of $A_{\tilde{s}}^i$ are less likely to explode.
Moreover, the staircase structure of $(A_s, B_s)$ is preserved by $(A_{\tilde{s}}, B_s)$, so the parametrization in Theorems~\ref{thm:observation-matrix-parametric-form} and~\ref{thm:parametric-Brunovsky-transformation} remain valid. 

Although our parametrization avoids explicit inversion of the controllability matrix, ill-conditioning of $T$ and $G$ can still degrade numerical stability. 
The parametric form provides a natural opportunity to mitigate this by directly optimizing their condition numbers.

The traditional $\kappa$-condition number, defined as the ratio between the largest and smallest singular values, is not an ideal objective function due to its non-smoothness.
The Frobenius-based condition number, though widely used, scales with the dimension and provides only a loose bound on the numerical errors in matrix inversion.
A recent alternative~\citep{jung2025omega} is the $\omega$-condition number, defined as the ratio between the arithmetic and geometric means of the singular values, which captures matrix conditioning while remaining differentiable. Their respective definitions are repeated below
\begin{equation}
    \kappa = \dfrac{\sigma_1}{\sigma_n}, \quad
    \omega = \dfrac{\frac{1}{n} \sum_{i=1}^n \sigma_i}{\Pi_{i=1}^n (\sigma_i)^\frac{1}{n}},
\end{equation}
where $\sigma_1 \geq \sigma_2 \geq \dots \geq \sigma_n > 0$ are the singular values. 

To avoid ill-conditioning, we minimize the condition numbers of the state transformation $T$ and the decoupling matrix $D$ in~\eqref{eq:def-parametric-Brunovsky-transformation}, since their inverses are required. 
The objective function is set as the logarithm of the product of their $\omega$-condition numbers for numerical performance.
The rank constraints in~\eqref{eq:rank-constraints} have measure zero and need no special treatment in the following unconstrained problem: 
\begin{align}\label{eq:opt-condition-number}
    \begin{split}
        \underset{S^r, S^f}{\text{min}} & \quad
        \log(\omega(T)) + \log(\omega(D)), \\
        \text{s.t. } & \quad T, D \text{ defined by~\eqref{eq:def-parametric-Brunovsky-transformation}}, \\
        & \quad S^r \text{ satisfy~\eqref{eq:rank-constraints}}.
    \end{split}
\end{align}

\section{Numerical experiments}\label{sec:num-exp}

We compare the proposed method (employing \texttt{fminunc} in \MATLAB to solve~\eqref{eq:opt-condition-number}) and Luenberger's method for constructing the Brunovsky transformation for $N=100$ randomly generated linear systems with $n=15$ states, $m = 4$ inputs, and controllability indices $\{\mu_1=\mu_2=5, \mu_3 =3, \mu_4 =2\}$. 
The road maps for the two methods are:
\begin{subequations}
\begin{align}
    &(A, B) \overset{\eqref{eq:trans-staircase}}{\rightarrow} (A_s, B_s) \overset{\eqref{eq:deadbeat-for-staircase}}{\rightarrow} (A_{\tilde{s}}, B_s) \overset{\eqref{eq:opt-condition-number}}{\rightarrow} (A_b, B_b) \tag{Proposed} \label{eq:proposed-method} \\
    &(A, B) \overset{\eqref{eq:trans-controllable-canonical}}{\rightarrow} (A_c, B_c) \overset{\text{last-rows}}{\underset{\text{elimination}}{\xrightarrow{\hspace{3em}}}} (A_b, B_b) \tag{Luenberger} \label{eq:Luenberger-method}
\end{align}
\end{subequations}

\begin{figure}
    \centering
    \includegraphics[width=1.0\linewidth]{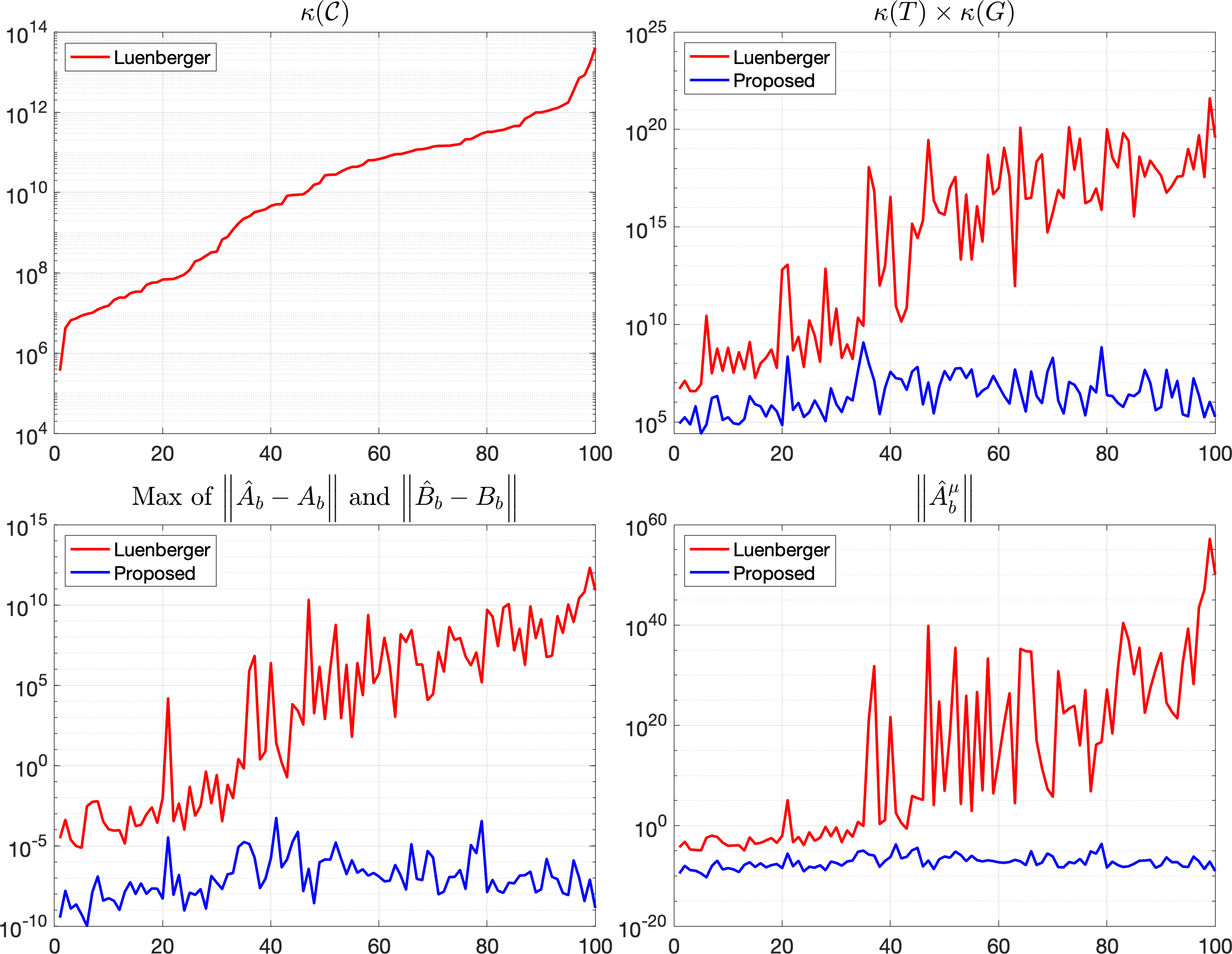}
    \caption{
    Red curve: \eqref {eq:Luenberger-method}. Blue curve: \eqref{eq:proposed-method}.
    Top left: condition number of the controllability matrix $\mathcal{C}$, characterizing how ill-conditioned the system is. 
    Top right: product of condition numbers of $T$ and $G$.
    Bottom left: the norm of the errors between the computed Brunovsky pair $(\hat{A}_b = T(A + BF)T\iv, \hat{B}_b=TBG)$ and the ground truth pair $(A_b, B_b)$. 
    Bottom right: the norm of $\hat{A}_b^\mu$. Ideally, $\norm{\hat{A}_b-A_b}$, $\norm{\hat{A}_b-A_b}$, and $\norm{\hat{A}_b^\mu}$ should all be zero, so a smaller value indicates a Brunovsky transformation with higher numerical reliability. 
    } 
    \label{fig:numerical-experiment}
\end{figure}

As shown in Fig.~\ref{fig:numerical-experiment}, as the linear systems become increasingly ill-conditioned, the classical Luenberger method produces huge errors in constructing the Brunovsky forms. The synthesized $\hat{A}_b$ is far from nilpotency. In contrast, the proposed method achieves errors between $10^{-5}$ and $10^{-10}$, guarantees nilpotency of $\hat{A}_b$, and maintains relatively well-conditioned $T$ and $G$, independent of the ill-conditioning of the original system.

We do not report the time spent by \texttt{fminunc} because the Brunovsky transformation is computed offline, prior to closed-loop operation. Therefore, the computational cost of solving~\eqref{eq:opt-condition-number} has a secondary impact, particularly for the low-dimensional systems typical of LTI MPC applications.

We compare against Luenberger because our method computes the full tuple $(T, F, G)$ rather than just the feedback gain. Existing software, such as SLICOT~\citep{benner1999slicot}, addresses pole placement alone and does not provide well-conditioned state and input transformations.

\section{Conclusion}

In this paper, we propose a numerically reliable method to construct Brunovsky transformations for general dense controllable linear systems in three steps: orthogonal similarity transformation to staircase form, feedback via a deadbeat control gain, and optimization over linearly parametrized transformations to minimize the condition number.
Numerical experiments demonstrate orders-of-magnitude improvements in construction error and conditioning compared with Luenberger’s classical method.
In future work, we plan to leverage the pseudo-convexity of $\kappa$ and $\omega$-condition numbers to establish convergence beyond local minima in~\eqref{eq:opt-condition-number}, validate the method on higher-dimensional systems, and prove Conjecture~\ref{conj:}.

\section*{DECLARATION OF GENERATIVE AI AND AI-ASSISTED TECHNOLOGIES IN THE WRITING PROCESS}
During the preparation of this work, the authors used ChatGPT and Claude to check grammar and improve writing. After using this tool/service, the authors reviewed and edited the content as needed and take full responsibility for the publication's content.

\bibliography{ifacconf}             % bib file to produce the bibliography
                                                     % with bibtex (preferred)

\appendix
\section{Proof of Lemma 2}

\begin{pf}
    Define $X_i = A_i A_{i-1} \dots A_2 \in \mathbb{R}^{r_i \times m}, i=p, p-1, \dots, 2$. All $X_i$ have full row rank by construction. 
    Since $A_1$ is invertible, it is equivalent to discuss $\tilde{D} = D A_1\iv$ instead. 
    \begin{equation}
        \tilde{D} = \begin{bmatrix}
            S_p A_p A_{p-1} \dots A_2 \\
            S_{p-1} A_{p-1} \dots A_2 \\
            \vdots \\
            S_2 A_2 \\
            S_1 
        \end{bmatrix}
    \end{equation}
    We first prove that if~\eqref{eq:rank-auxiliary-lemma-constraints} is satisfied, then $\tilde{D}$ is invertible. 
    $S_p$ is invertible, so the top block $S_p X_p$ has rank $r_p$. 
    $\begin{bsmallmatrix}
        A_p \\
        S_{p-1}
    \end{bsmallmatrix}$ is invertible, so the row spaces of $A_p$ and $S_{p-1}$ are independent. So the top two block of $\tilde{D}$, $\begin{bsmallmatrix}
        S_p A_p \\
        S_{p-1}
    \end{bsmallmatrix}X_{p-1}$, has rank $r_{p-1}$. 
    The logic continues for the top three blocks 
    $\begin{bsmallmatrix}
        S_p A_p A_{p-1} \\
        S_{p-1} A_{p-1} \\
        S_{p-2}
    \end{bsmallmatrix}X_{p-2}$ having rank $r_{p-2}$, until the top $p$ block rows having rank $r_0 = m$, which means $\tilde{D}$ is invertible. \\
    We then prove that if $\tilde{D}$ is invertible, \eqref{eq:rank-auxiliary-lemma-constraints} must be satisfied. 
    The top block $S_p X_p$ must have row rank $r_p$, so $S_p$ is invertible. 
    The top two blocks $\begin{bsmallmatrix}
        S_p A_p \\
        S_{p-1}
    \end{bsmallmatrix}X_{p-1}$ must have row rank $r_{p-1}$, so $\begin{bsmallmatrix}
        S_p A_p \\
        S_{p-1}
    \end{bsmallmatrix}$ is invertible $\Rightarrow \begin{bsmallmatrix}
        A_p \\
        S_{p-1}
    \end{bsmallmatrix}$ is invertible because $S_p$ does not change the row space of $A_p$. 
    The top three blocks $\begin{bsmallmatrix}
        S_p A_p A_{p-1} \\
        S_{p-1} A_{p-1} \\
        S_{p-2}
    \end{bsmallmatrix}X_{p-2}$ must have row rank $r_{p-2} \Rightarrow S_{p-2}$ must span the $k_{p-2}$ dimensional row spaces that is uncovered by $A_{p-1} \Rightarrow \begin{bsmallmatrix}
        A_{p-1} \\
        S_{p-2}
    \end{bsmallmatrix}$ is invertible. 
    The same reasoning extends to the final case: the top $p$ blocks, $\tilde{D}$, are of rank $r_0=m$, so $\begin{bsmallmatrix}
        A_2 \\
        S_1
    \end{bsmallmatrix}$ is invertible. 
    \hfill $\blacksquare$
\end{pf}

\end{document}